\input amstex
\magnification=1200
\documentstyle{amsppt}
\NoRunningHeads
\NoBlackBoxes
\define\CA{\Cal A}
\define\CX{\Cal X}
\define\HCA{\widehat{\CA}}
\define\CU{\Cal U}
\define\CH{\Cal H}
\define\CDb{\Cal D}
\define\CG{\Cal G}
\define\CK{\Cal K}
\define\SU{\operatorname{SU}}
\define\ga{\alpha}
\define\sL{\operatorname{\frak s\frak l}}
\define\sO{\operatorname{\frak s\frak o}}
\define\Fg{\frak g}
\define\Fh{\frak h}
\define\Fb{\frak b}
\define\Fq{\frak q}
\define\BE{\bold E}
\define\Aut{\operatorname{Aut}}

\define\osc{\operatorname{\frak o\frak s\frak c}}
\document\qquad\qquad\qquad\qquad\qquad\qquad\qquad\qquad\qquad\qquad
$\boxed{\boxed{\aligned
&\text{\eightpoint"Thalassa Aitheria" Reports}\\
&\text{\eightpoint RCMPI-96/06$^+$ (November 1996)}\endaligned}}$\newline
\ \newline
\ \newline
\topmatter
\title Topics in hidden symmetries. III.
\endtitle
\author Denis V. Juriev
\endauthor
\affil\eightpoint\rm
"Thalassa Aitheria" Research Center for Mathematical Physics and
Informatics,\linebreak
ul.Miklukho-Maklaya 20-180, Moscow 117437 Russia.\linebreak
E-mail: denis\@juriev.msk.ru
\endaffil
\date E-version: q-alg/yymmxxx
\enddate
\abstract This short paper being devoted to some aspects of the inverse
problem of the representation theory treats several themes, which have their
origins in the researches of F.A.Berezin, D.P.Zhelobenko, V.P.Maslov and his
group, in context of the author's approach to the setting free of hidden
symmetries, which is based on the noncommutative geometry.

Relations of the discussed topics to the general ideologies of the translator
algebras in representation theory and of the dynamical symmetries or the
spectrum--generating algebras in mathematical physics are pointed out.
A possibility of quantum field generalizations of the constructions in terms
of vertex (QFT--operator) and Zamolodchikov algebras is mentioned.
\endabstract
\endtopmatter
This paper being a continuation of the previous two parts [1,2] is a collection
of examples as illustrating the general ideology presented in the review [3]
as explicating its new features. Similar to the part II all examples emphasize
the interrelations between the material of the previous papers on the setting
free of hidden symmetries [4,1,3:\S1], the quantization of constants [3:\S3]
and the researches of D.P.Zhelobenko on the generalized Mickelsson algebras
[5,6] and related objects of the representation theory (i.e. the algebraic
extremal projector, see also [7]), so the material may be considered also as
the author's variations on themes of the book [6]. Also the links of the
inverse problem of representation theory with noncommutative geometry,
which were pointed out in [3], are interpreted in a new way. Briefly, the
inverse problem is considered in the context of approach to noncommutative
geometries of A.Connes [8] and Yu.Manin [9], which naturally appears in the
representation theory [6] (the noncommutative quadratic algebras over the
functional rings and fields in sense of D.P.Zhelobenko), such approach is also
related to the theory of nonlinear Poisson brackets and asymptotic quantization
[10] and to the theory of classical and quantum dynamical systems (see f.e.
[11]).

The concrete lines of topics have their origins in the author's papers
[12-14], where some objects, which gave start to the constructions below,
appeared.

\head 1. Topic Six: A variation on two themes
\endhead

This topic is devoted to a variation on two themes, the first of them is
taken from the book [6] of D.P.Zhelobenko and the second one belongs to the
author [3:\S3]. The variation is also penetrated by tunes from the book [10]
(as well as from the elder review [15]) and from the articles [16,17], arranged
in spirit of the general author's ideas on the setting free of hidden
symmetries [1-3].

\definition{Definition 1 (cf.[6])} Let $\CX$ be an associative algebra, then
the associative algebra $\CA$ is called {\it an algebra over $\CX$ (in
sense of D.P.Zhelobenko)\/} iff $\CA$ is the left $\CX$--module, whose
structure is compatible with multiplication in $\CA$ (i.e.
$\forall x\in\CX$, $\forall a,b\in\CA$ $x(ab)=(xa)b$) and such that
there exist a linear mapping $\mu:\CA\otimes\CX\mapsto\CX\otimes\CA$,
$\mu(a,x)=\sum_{\ga}x^{(\ga)}\otimes a^{(\ga)}$ so that
$\forall a,b\in\CA$, $\forall x\in\CX$ $a(xb)=\sum_{\ga}x^{(\ga)}(a^{(\ga)}b)$.
\enddefinition

\remark{Remark 1} One may consider the operator product expansions, R--matrices,
QFT--operator algebras (vertex algebras) and Zamolodchikov algebras over
functional commutative algebras and fields (in sense of D.P.Zhelobenko) and,
even, over their noncommutative analogs.

Such Zamolodchikov algebras over the commutative functional algebra (generated
by the momenta) systematically appear in the R--matix constructions of the
quantum Liouville field theory (works of J.-L.Gervais and collaborators on the
quantum group aspects of the theory).

The reader should keep in mind that though we have deal with simple examples,
which are easily computable, it is very desirable to generalize all
constructions of this topics (as well as of certain others) to the
parametric case.
\endremark

\definition{Definition 2}

{\bf A.} Let $\CA$ and $\CX$ be two arbitrary Poincar\`e--Birkhoff--Witt
algebras (PBW--algebra [10:App.2], see also [15]) with basises $e_i$
($i=1,\ldots n$), $x_j$ ($j=1,\ldots m$) and commutation relations
$[e_i,e_j]=f_{ij}(e_1,\ldots e_n)$, $[x_i,x_j]=g_{ij}(x_1,\ldots x_m)$.
The PBW--algebra $\HCA$ with basis $a_i$ ($i=1,\ldots n$), $\eta_j$
($j=1,\ldots m$) is called {\it the partial linearization of $\CA$ over
$\CX$\/} iff:
\roster
\item there exist an epimorphism $\pi:\HCA\mapsto\CA$ such that
$\pi(a_i)=e_i$ and monomorphism $\iota:\CX\mapsto\HCA$ such that
$\iota(x_j)=\eta_j$;
\item the monomorphism $\eta$ supply $\HCA$ by the structure of the algebra
over $\CX$ (in sense of D.P.Zhelobenko);
\item $[a_i,a_j]=h_{ij}^k(\eta_1,\ldots\eta_m)a_k$ (i.e. $\HCA$ is a Lie
algebra over $\CX$ (in sense of D.P.Zhelobenko).
\endroster

{\bf B.} Let $\CA$ and $\CX$ be the same data. The PBW--algebra $\HCA$ with
basis $a_i$ ($i=1,\ldots n$), $\eta_j$ ($j=1,\ldots m$) is called {\it the
generalized partial linearization of $\CA$ over $\CX$\/} iff:
\roster
\item there exist an epimorphism $\pi:\HCA\mapsto\CA$ such that
$\pi(a_i)=e_i$ and monomorphism $\iota:\CX\mapsto\HCA$ such that
$\iota(x_j)=\eta_j$;
\item the monomorphism $\eta$ supply $\HCA$ by the structure of the algebra
over $\CX$ (in sense of D.P.Zhelobenko);
\item the algebra $\HCA$ is the algebra with quadratic (perphaps,
non-homogeneous) relations over $\CX$ (in sense of D.P.Zhelobenko).
\endroster

{\bf C.} Let $\CA$ and $\CX$ be the same data, $\CK$ be an extension of $\CX$.
The algebra $\CA$ is called $(\CX,\CK)$--regular iff there exists its
generalized partial linearization $\HCA$ over $\CK$ such that for each its
generalized partial linearization $\HCA'$ over $\CX$ there exists a
monomorphism $\vartheta:\CX\mapsto\CK$ and a monomorphism of the algebra
$\HCA'$ over $\CX$ into the algebra $\HCA$ over $\CK$, which is considered as
an algebra over $\vartheta(\CX)$. The algebra $\CA$ over the Ore algebra
$\CX$ is called $\CX$--regular iff it is $(\CX,\CDb(\CX))$--regular, where
$\CDb(\CX)$ is the algebra of fractions for $\CX$.
\enddefinition

\example{Example 1} Let $\CA$ be the nonlinear $\sL_2$ [18] with the basis
$e_i$ ($i=-1,0,1$) and the commutation relations $[e_{\pm1},e_0]=\pm e_{\pm1}$,
$[e_1,e_{-1}]=h(e_0)$, where $h(t)=th_0(t)$. Let $\CX$ be an algebra generated
by a variable $x$; $\HCA$, the partial linearization of $\CA$ over $\CX$, is
the algebra generated by $a_i$ ($i=-1,0,1$), $\eta$ with the commutation
relations $[a_{\pm1},a_0]=\pm a_{\pm1}$, $[a_1,a_{-1}]=h_0(\eta)a_0$,
$[a_i,\eta]=ia_i$ and such that $\pi(a_i)=e_i$, $\pi(\eta)=e_0$,
$\iota(x)=\eta$.
\endexample

\example{Example 2} Let $\CA=\CU_q(\sL_2)$ [17;6:Ch.9] with the basis $e_i$
($i=-1,0,1$) and the commutation relations $[e_{\pm1},e_0]=\pm e_{\pm1}$,
$[e_1,e_{-1}]=\frac{q^{e_0}-q^{-e_0}}{q-q^{-1}}$. Let $\CX$ be an algebra
generated by two commuting variables $x_{\pm}$; $\HCA$, the partial
linearization of $\CA$ over $\CX$, is the algebra generated by $a_i$
($i=-1,0,1$), $\eta_{\pm}$ with commutation relations $[a_{\pm1},a_0]=\pm
a_{\pm1}$, $[a_1,a_{-1}]=\eta_+-\eta_-$, $a_i\eta_{\pm}=q^{\pm i}\eta_{\pm}a_i$,
$[\eta_+,\eta_-]=0$ and such that $\pi(a_i)=e_i$,
$\pi(\eta_{\pm})=\frac{q^{\pm e_0}}{q-q^{-1}}$, $\iota(x_{\pm})=\eta_{\pm}$.
\endexample

\definition{Definition 3 (cf.[3:\S3])}

{\bf A.} Let $\Fg$ be a Lie algebra with the fixed basis $e_i$ ($i=1,\ldots n$)
and commutation relations $[e_i,e_j]=c^k_{ij}e_k$ and $\CX$ be a PBW--algebra
with the basis $x_j$ ($j=1,\ldots m$) and commutation relations
$[x_i,x_j]=h_{ij}(x_1,\ldots x_m)$. A PBW--algebra $\CA$ with the basis
$a_i$ ($i=1,\ldots n$) $\eta_j$ ($j=1,\ldots m$) is the result of {\it the
quantization of constants\/} in $\Fg$ with {\it the algebra of quantized
constants\/} $\CX$ iff
\roster
\item there exist a monomorphism $\iota:\CX\mapsto\CA$, $\iota(x_i)=\eta_i$;
\item $[a_i,a_j]=f_{ij}^k(\eta_1,\ldots \eta_m)$ and that
$(c_{ij}^k=0)\Rightarrow(f_{ij}^k\equiv0)$.
\endroster
The mapping $\Fq:\Fg\mapsto\CA$, $\Fq(e_i)=a_i$ is called {\it the quantization
of constants}.

{\bf B.} The quantization of constants is called {\it quasilinear\/} iff
the monomorphism $\iota$ supplies $\CA$ by the structure of algebra over
$\CX$ (in sense of D.P.Zhelobenko).

{\bf C.} Let $\{\Fg_{\ga}\}$ be a set of subalgebras in the Lie algebra
$\Fg$. The quantization of constants is called $\{\Fg_{\ga}\}$--{\it preserving\/}
iff $\Fq$ realizes a monomorphism of each $\Fg_{\ga}$ into the commutator
algebra $\CA_{[\cdot,\cdot]}$ of the algebra $\CA$.
\enddefinition

Many generalized partial linearizations may be regarded as results of
the quantization of constants in Lie algebras.

\example{Example 3} The algebras of examples 1 and 2 are the results of
the quantization of constants in the Lie algebra $\sL_2$. They are quasilinear
and $\Fb_{\pm}$--preserving, where $\Fb_{\pm}$ are Borel subalgebras of
$\sL_2$ generated by $e_0$ and $e_{\pm}$, respectively.
\endexample

\example{Example 4 {\rm [3:\S3.1]}} The Sklyanin algebra is the result of
quantization of constants in $\sO_3$, which is not quasilinear.
\endexample

\example{Example 5: the Lobachevski{\v\i} algebra
(cf.[12:\S3.2;13:\S2.2;14:App.A.3])}\linebreak
The relations between the author's ideas on the setting free of hidden
symmetries and the Lobachevski{\v\i} algebra (the Berezin quantization of
the Lobackevski{\v\i} plane [16], recently rediscovered by S.Klimek and
A.Lesniewski [19]) are discussed below.

In the Poincare realization of the Lobachevskii plane (the realization in the
unit disk) the Lobachevskii metric may be written as
$w=q_R^{-1}\,dzd\bar{z}/(1-|z|^2)^2$;
one can construct the $C^*$--algebra (Lobachevski{\v\i} algebra), which may be
considered as a quantization of such metric, namely, let us consider
two variables $t$ and $t^*$, which obey the following commutation relations:
$[tt^*,t^*t]=0$, $[t,t^*]=q_R(1-tt^*)(1-t^*t)$ (or in an equivalent form
$[ss^*,s^*s]=0$, $[s,s^*]=(1-q_Rss^*)(1-q_Rs^*s)$, where $s=(q_R)^{-1/2}t$).

One may realize such variables by tensor operators in the Verma module over
$\sL_2$ of the weight $h=\frac{q_R^{-1}+1}2$; if such Verma module is realized
in polynomials of one complex variable $z$ and the action of $\sL_2$ has the
form $L_{-1}=z$, $L_0=z\partial_z+h$, $L_1=z(\partial_z)^2+2h\partial_z$, then
the variables $t$ and $t^*$ are represented by tensor operators $D=\partial_z$
and $F=z/(z\partial_z+2h)$, where $[L_i,D]=-D^{i+1}$, $[L_i,F]=F^{i-1}$. These
operators are bounded if $q_R>0$ and therefore one can construct a Banach
algebra generated by them and obeying the prescribed commutation relations;
the structure of $C^*$--algebra is rather obvious: an involution $*$ is
defined on generators in a natural way, because the corresponding tensor
operators are conjugate to each other.

Let us consider a variable $x$ and the algebra $\CX=\Bbb C[x]$ of polynomials
of $x$. A partial linearization of the Lobachevski{\v\i} $C^*$--algebra over
$\CX$ may be constructed; it is generated by the elements $\eta$, $\tau$,
$\tau^*$, $\iota(x)=\eta$, $\pi(\tau)=t$, $\pi(\tau^*)=t^*$, $\pi(\eta)=
q_R(1-tt^*)$ and the commutation relations have the form
$$[\tau,\tau^*]=q_R^{-1}\eta^2(1-\eta)^{-1},\quad
\eta\tau=\tau\eta(1-\eta)^{-1},\quad
\eta\tau^*=\tau^*\eta(1+\eta)^{-1}.$$
The obtained algebra is result of the certain quantization of a constant in the
Heisenberg algebra.

It is possible to construct also a family of generalized partial
linearizations of the Lo\-ba\-chev\-ski{\v\i} algebra. One of them is generated
by the elements $\eta$, $\tau$, $\tau^*$, $\pi(\tau)=t$, $\pi(\tau^*)=t^*$,
$\pi(\eta)=q_R(1-tt^*)$ and the commutation relations have the form
$$\tau\tau^*-(1-\eta)\tau^*\tau=\eta,\quad
\eta\tau=\tau\eta(1-\eta)^{-1},\quad
\eta\tau^*=\tau^*\eta(1+\eta)^{-1}.$$
If one put $\xi=\eta^{-1}$ then the commutation relations will be rewritten as
$$\xi\tau\tau^*-(\xi-1)\tau^*\tau=1,\quad [\xi,\tau]=-\tau,\quad
[\xi,\tau^*]=\tau^*.$$
\endexample

\proclaim{Theorem 1} The Lobachevski{\v\i} algebra is $\Bbb C[x]$--regular.
\endproclaim

The corresponding universal algebra $\HCA$ over the algebra of rational
functions $\Bbb C(x)=\CDb(\Bbb C[x])$ may be constructed from the oscillator
Lie algebra. The oscillator Lie algebra $\osc$ is generated by the elements
$p$, $q$, $r$, $\epsilon$, where $p$, $q$ and $r$ form the Heisenberg
subalgebra ($[p,q]=r$, $[r,p]=[r,q]=0$) and $[\epsilon,q]=q$, $[\epsilon,p]=-p$,
$[\epsilon,r]=0$ so that $\HCA=\CU'(\osc)=\CU(\osc)\otimes_{\Bbb
C[\epsilon]}\Bbb C(\epsilon)$.

\head 2. Topic Seven: A theme more\endhead

In this topic we shall incorporate one theme more into our picture. The main
objects of this theme are generalized Mickelsson and Zhelobenko algebras [5,6]
as well as the extremal projectors [7], so the exposition is plunged into
the ideologies of translator algebras in the representation theory (see
f.e.[20,6]) and of dynamical symmetries or the spectrum--generating
algebras in mathematical physics [21].

\definition{Definition 4 (a preliminary simplified version)} Let
$\CX=D(\Bbbk^n)$ ($\Bbbk$ is the basic field, $\Bbbk^n$ is the abelian Lie
algebra and $D$ denotes the algebra of fractions of the Lie algebra [20]),
$\HCA$ be an algebra over $\CX$ (in sense of D.P.Zhelobenko), $\Fg$ be a
reductive Lie algebra, $\operatorname{rank}\Fg=n$. The envelopping algebra
$\BE(\HCA)$ of $\Fg$ is called {\it the extremal extension of\/} $\HCA$ iff the
pair $(\BE(\HCA),\Fg)$ obey the Zhelobenko conditions [6:\S7.3] and $\HCA$ is
isomorphic to the Zhelobenko algebra $Z(\BE(\HCA),\Fg)$ [5,6], here $\CX$ is
identified with $D(\Fh)$ ($\Fh$ is the Cartan subalgebra of $\Fg$).
\enddefinition

\remark{Remark 2} This definition has a preliminary character, its correct
formulation should be based on the final complete axiomatics of the convenient
class of contragredient algebras [22;6:App.A], whose objects should be used
in the definition instead of $\Fg$. One of the general versions is formulated
below.
\endremark

\definition{Definition 4 (a currently used general version)} Let $\CH$ be
a regular commutative algebra (i.e. an algebra with the set of characters
being total in it [6:App.A]) without zero divisors, $Q$ be a commutative
subgroup of $\Aut(\CH)$ and $\CX$ be a commutative extension of $\CH$ such that
(1) $\CX$ has no zero divisors, (2) $\CDb(\CH)\subseteq\CX$ (here $\CDb$
denotes the algebra of fractions of an associative algebra),
(3) $\Aut(\CH)\subseteq\Aut(\CX)$. Let also $\HCA$ be an algebra over $\CX$ and
$\CG$ be a regular contragredient (and, hence, weakly triangular) algebra of
the Cartan type [6:App.A] with the Cartan subalgebra $\CH$. The finitely
generated unital algebra $\BE(\HCA)$ over $\CG$ (in sense of D.P.Zhelobenko)
with quadratic relations over $\CG$ is called {\it the extremal extension of\/}
$\HCA$ iff $\HCA$ is isomorphic to the Zhelobenko algebra
$Z(\BE(\HCA)^{\sigma}(\CX),\CG)$, where $\BE(\HCA)^{\sigma}(\CX)$ is the second
canonical extension [6:App.A] of the algebra $\BE(\HCA)$.
\enddefinition

Below we shall be interested in the extremal extensions $\BE(\HCA')$ of
algebras $\HCA'=\HCA\otimes_{\CDb(\CX)}\CX$, where $\HCA$ are (generalized)
partial linearizations of algebras $\CA$ over $\CX$ (briefly, $\BE(\HCA')$
will be called the extremal extension of the algebra $\CA$ and be denoted by
$\BE(\CA)$).

\proclaim{Theorem 2 {\it (Example 6: $\CU(\sL_2\oplus\sL_2)$ as an extremal
extension of the Lobachevski{\v\i} algebra)}} The universal envelopping
algebra $\CU(\sL_2\oplus\sL_2)$ is an extremal extension of the Lobachevski{\v\i}
algebra.
\endproclaim

\demo{Proof} The Zhelobenko algebra $Z(\CU(\sL_2\oplus\sL_2),\sL_2)$
was constructed in [6:\S2.4.6]. It is generated by three variables $s_0$,
$s_-$ and $s_+$ ($s_0$ belongs to the center). After the factorization over
the central commutative subalgebra $\Bbb C[s_0]$ it is isomorphic to
$\CU'(\osc)$, which is the universal generalized partial linearization of
the Lobachevski{\v\i} algebra (see Theorem 1) \qed
\enddemo

\remark{Remark 3} The embedding of the generalized partial linearizations
of the Lo\-ba\-chev\-s\-ki{\v\i} algebra into the factorized Zhelobenko algebra
$Z(\CU(\sL_2\oplus\sL_2),\sL_2)$ defines the central extension of the
Lobachevski{\v\i} algebra, whose generalized partial linearizations are
embed into the Zhelobenko algbera itself.
\endremark

\Refs
\roster
\item" [1]" Juriev D., Topics in hidden symmetries. I. E-print: hep-th/9405050.
\item" [2]" Juriev D., Topics in hidden symmetries. II. Report RCMPI-96/06
(October 1996) [e-version: q-alg/9610026 (1996)].
\item" [3]" Juriev D., An excursus into the inverse problem of representation
theory. Report RCMPI-95/04 (1995) [e-version: mp\_arc/96-477].
\item" [4]" Juriev D., Setting hidden symmetries free by the noncommutative
Veronese mapping. J.Math.Phys. 35(9) (1994) 5021-5024.
\item" [5]" Zhelobenko D.P., Extremal projectors and generalized Mickelsson
algebras over reductive Lie algebras [in Russian]. Izvestiya AN SSSR.
Ser.matem. 52(4) (1988) 758-773.
\item" [6]" Zhelobenko D.P., Representations of the reductive Lie algebras.
Moscow, Nauka, 1994.
\item" [7]" Asherova R.M., Smirnov Yu.F., Tolsto\v\i\ V.N., Projector operators
for simple Lie algebras [in Russian]. Teor.Matem.Fiz. 8(2) (1971) 255-271;
A description of a class of projector operators for semisimple complex Lie
algebras [in Russian]. Matem.Zametki 26(1) (1979) 15-26.
\item" [8]" Connes A., Introduction \`a la g\'eometrie noncommutative.
InterEditions, Paris, 1990.
\item" [9]" Manin Yu.I., Topics in noncommutative geometry, Princeton Univ.
Press, Princeton, 1991.
\item"[10]" Karasev M.V., Maslov V.P., Nonlinear Poisson brackets. Geometry and
quantization. Amer.Math.Soc., Providence, RI, 1993.
\item"[11]" Latushkin Yu.D., Stepin A.M., Weighted shift operators and linear
extensions of dynamical systems [in Russian]. Uspekhi Matem.Nauk. 46 (1991)
85-143.
\item"[12]" Juriev D.V., Quantum projective field theory: quantum--field
analogs of the Euler--Arnold equations in projective $G$--hypermultiplets [in
Russian]. Teor.Matem.Fiz. 98(2) (1994) 220-240 [English transl.:
Theor.Math.Phys. 98 (1994) 147-161].
\item"[13]" Juriev D.V., Complex projective geometry and quantum projective
field theory [in Russian]. Teor.Matem.Fiz. 101(3) (1994) 331-348 [English
transl.: Theor.Math.Phys. 101 (1994) 1387-1403].
\item"[14]" Juriev D.V., Belavkin--Kolokoltsov watch--dog effects in
interactively controlled stochastic dynamical videosystems [in Russian].
Teor.Matem.Fiz. 106(2) (1996) 333-352 [English transl.: Theor.Math.Phys.
106 (1996) 000-000].
\item"[15]" Karasev M.V., Maslov V.P., Naza{\v\i}kinski{\v\i} V.E., Algebras
with general commutation relations [in Russian]. Current Probl. Math.,
Modern Achievements. V.13. Moscow, VINITI, 1979.
\item"[16]" Berezin F.A., Quantization in complex symmetric spaces [in
Russian]. Izvestiya AN SSSR. Ser.matem. 39(2) (1975) 363-402.
\item"[17]" Reshetikhin N.Yu., Takhtadzhyan L.A., Faddeev L.D.,
Quantization of Lie algebras and Lie groups. Algebra i anal. 1(2) (1989)
178-206 [English transl.: St.Petersburg Math.J. 1 (1990) 193-225].
\item"[18]" M. Ro\v cek, Representation theory of the nonlinear $\SU(2)$
algebra. Phys. Lett.B 255 (1991) 554-557.
\item"[19]" Klimek S., Lesniewski A., Commun.Math.Phys. 146 (1992) 103-122.
\item"[20]" Dixmier J., Alg\`ebres enveloppantes. Paris, Villars, 1973.
\item"[21]" Barut A., Raczka R., Theory of group representations and
applications. V.II., PWN--Polish Sci. Publ., Warszawa, 1977;\newline
Symmetries in science. VII. Spectrum--generating algebras and dynamic
symmetries in physics. Eds.B.Gruber, T.Otsake, Plenum Publ., 1995.
\item"[22]" Zhelobenko D.P., Contragredient algebras. J.Group Theory Phys.,
1(1) (1993) 201-233.
\endroster
\endRefs
\enddocument